# TIME AVERAGES, RECURRENCE AND TRANSIENCE IN THE STOCHASTIC REPLICATOR DYNAMICS

By Josef Hofbauer and Lorens A. Imhof

*University of Vienna and Bonn University*

We investigate the long-run behavior of a stochastic replicator process, which describes game dynamics for a symmetric two-player game under aggregate shocks. We establish an averaging principle that relates time averages of the process and Nash equilibria of a suitably modified game. Furthermore, a sufficient condition for transience is given in terms of mixed equilibria and definiteness of the payoff matrix. We also present necessary and sufficient conditions for stochastic stability of pure equilibria.

**1. Introduction.** The present paper deals with a stochastic variant of the continuous-time replicator dynamics. We begin with a brief review of the deterministic model. For a comprehensive discussion, see Hofbauer and Sigmund (1998), Weibull (1995) and the recent surveys by Hofbauer and Sigmund (2003), Nowak and Sigmund (2004) and Sandholm (2007). Consider a symmetric two-player game with $n$ pure strategies, $1, \ldots, n$, and $n \times n$ payoff matrix $A = (a_{ij})$. Thus both players have the same set of strategies and for either player, $a_{ij}$ is the payoff from using strategy $i$ if the opponent uses strategy $j$. There is no symmetry or skew-symmetry assumption on $A$. The replicator dynamics describes how the proportions of strategies in a population evolve. Consider a large population where every member is programmed to play one pure strategy. Let $\zeta_i(t)$ denote the size of the subpopulation of $i$-players at time $t$, and let $\xi_i(t) = \zeta_i(t)/[\zeta_1(t) + \cdots + \zeta_n(t)]$ denote its proportion. If the population is in state $\boldsymbol{\xi}(t) = (\xi_1(t), \ldots, \xi_n(t))^T$, then, under random matching, $\{A\boldsymbol{\xi}(t)\}_i$ is the expected payoff to individuals playing $i$. Suppose that this payoff represents the per capita growth rate in the $i$th subpopulation. Thus

$$\dot{\zeta}_i = \zeta_i \{A\boldsymbol{\xi}\}_i, \qquad i = 1, \ldots, n. \tag{1.1}$$

This yields the deterministic replicator dynamics of Taylor and Jonker (1978)

$$(1.2) \qquad \dot{\xi}_i = \xi_i[\{A\boldsymbol{\xi}\}_i - \boldsymbol{\xi}^T A \boldsymbol{\xi}], \qquad i=1,\ldots,n.$$

The observation that biological processes modeled by replicator dynamics are inherently stochastic in nature led Foster and Young (1990) to introduce a replicator model based on a stochastic differential equation, and it turned out that even small stochastic effects can qualitatively change the asymptotic behavior. The present paper studies a stochastic model of Fudenberg and Harris (1992) which is similar to that of Foster and Young but has a more natural boundary behavior. Following Fudenberg and Harris, we introduce random perturbations to the payoffs modeled by independent Gaussian white noises with intensities $\sigma_1^2, \ldots, \sigma_n^2$. Instead of (1.1) we consider

$$(1.3) \qquad dZ_i(t) = Z_i(t)[\{AX(t)\}_i\, dt + \sigma_i\, dW_i(t)], \qquad i=1,\ldots,n,$$

where

$$X = (X_1, \ldots, X_n)^T = \frac{1}{Z_1 + \cdots + Z_n}(Z_1, \ldots, Z_n)^T$$

and $(W_1, \ldots, W_n)^T = W$ is an $n$-dimensional Brownian motion. The evolution of the population state $X(t)$ is then given by the stochastic replicator dynamics

$$(1.4) \qquad dX(t) = \mathbf{b}(X(t))\, dt + C(X(t))\, dW(t), \qquad X(0) = \mathbf{x}_0,$$

where

$$\mathbf{b}(\mathbf{x}) = [\mathrm{diag}(x_1, \ldots, x_n) - \mathbf{x}\mathbf{x}^T][A - \mathrm{diag}(\sigma_1^2, \ldots, \sigma_n^2)]\mathbf{x}$$

and

$$C(\mathbf{x}) = [\mathrm{diag}(x_1, \ldots, x_n) - \mathbf{x}\mathbf{x}^T]\mathrm{diag}(\sigma_1, \ldots, \sigma_n)$$

for $\mathbf{x} \in \Delta = \{\mathbf{y} \in [0,1]^n : y_1 + \cdots + y_n = 1\}$, and $\mathbf{x}_0$ is an initial value in int $\Delta$.

Fudenberg and Harris give a complete analysis of the asymptotic behavior for games with two pure strategies. Further papers that study the case $n=2$ include Saito (1997), Amir and Berninghaus (1998), Corradi and Sarin (2000) and Beggs (2002). Most of the analysis for the case $n=2$ relies on tools specific to one-dimensional diffusions. Few papers treat the general case $n \geq 2$. Cabrales (2000) considers extinction of dominated strategies, Imhof (2005a) examines the long-run behavior in the presence of evolutionarily stable strategies and Khasminskii and Potsepun (2006) analyze a related stochastic replicator model with Stratonovich type random perturbations. Benaïm, Hofbauer and Sandholm (2008) study the relation between recurrence and permanence. Imhof (2008) analyzes stochastic dynamics for games that describe multiple-trial conflicts. Discrete stochastic replicator processes have recently been investigated by Schreiber (2001) and Benaïm, Schreiber and Tarrès (2004).

The aim of the present paper is to provide further insight into the long-run behavior of the stochastic process $\{X(t)\}$ and in particular of its time averages $T^{-1} \int_0^T X(t)\,dt$ in the general case $n \geq 2$. It turns out that a crucial role is played by the modified payoff matrix

(1.5) $$\tilde{A} = (\tilde{a}_{ij})_{i,j=1}^n, \qquad \tilde{a}_{ij} = a_{ij} - \tfrac{1}{2}\sigma_i^2.$$

In contrast to the deterministic solution $\boldsymbol{\xi}(t)$, the stochastic process $\{X(t)\}$ cannot converge to an interior point of $\Delta$. In fact, the only points to which $\{X(t)\}$ can converge with positive probability are the vertices of $\Delta$, which correspond to populations consisting of one single type of players. This is one reason for our interest in the analysis of the time averages $T^{-1}\int_0^T X(t)\,dt$, even in cases where (1.2) has a globally asymptotically stable equilibrium. In Section 2 we study limit points of the time averages and relate them to Nash equilibria and correlated equilibria of the modified game $\tilde{A}$.

Section 3 deals with recurrence. As an immediate consequence of the results in Section 2, we obtain the following averaging principle: if $\{X(t)\}$ is positive recurrent, then $\tilde{A}$ has a unique interior Nash equilibrium and the time averages converge almost surely to this equilibrium. Furthermore, we consider a certain class of games that includes all zero-sum games, and give an explicit expression of the invariant densities for the transition probability functions of the corresponding replicator processes.

In Section 4 we derive an exclusion principle: if there does not exist a strategy against which every opponent obtains the same payoff under $\tilde{A}$, then $\{X(t)\}$ is transient. It is also shown that the process is transient if there does exist a mixed strategy with this property and $\tilde{A}$ is conditionally positive definite, that is, $\mathbf{y}^T \tilde{A} \mathbf{y} > 0$ for every $\mathbf{y} \neq \mathbf{0}$ with $y_1 + \cdots + y_n = 0$. Moreover, we investigate the relation between pure Nash equilibria of $\tilde{A}$ and stochastic stability of the corresponding vertices of $\Delta$. Thus we obtain three sufficient conditions for transience of the replicator process. Section 5 relates our results to some of the central results for the deterministic replicator dynamics, including the folk theorem of evolutionary game theory.

If we introduce random perturbations into (1.1) using Stratonovich integrals instead of Itô integrals, we obtain

(1.6) $\quad dZ_i^{(s)}(t) = Z_i^{(s)}(t)\{AX^{(s)}(t)\}_i\,dt + \sigma_i Z_i^{(s)}(t) \circ dW_i(t), \qquad i = 1, \ldots, n,$

where $X^{(s)} = (Z_1^{(s)} + \cdots + Z_n^{(s)})^{-1} Z^{(s)}$. See Turelli (1977) for a comparison of Itô and Stratonovich equations from a biological point of view. The Stratonovich equation (1.6) is equivalent to the Itô equation

$$dZ_i^{(s)}(t) = Z_i^{(s)}(t)[\{AX^{(s)}(t)\}_i + \tfrac{1}{2}\sigma_i^2]\,dt + \sigma_i Z_i^{(s)}(t)\,dW_i(t), \qquad i = 1,\ldots,n.$$

Thus $\{X^{(s)}(t)\}$ can be regarded as a solution to the Itô stochastic replicator equation of Fudenberg and Harris for the payoff matrix $A^{(s)}$ with entries $a_{ij}^{(s)} = a_{ij} + \frac{1}{2}\sigma_i^2$. Consequently, all the results in the present paper carry over to the Stratonovich solution $\{X^{(s)}(t)\}$ when $A$ is replaced by $A^{(s)}$. For instance, Theorem 3.1(a) yields that if $\{X^{(s)}(t)\}$ is positive recurrent, then the time averages $T^{-1} \int_0^T X^{(s)}(t)\,dt$ converge to an interior Nash equilibrium of $A$. Note that modifying $A^{(s)}$ as in (1.5) leads to the matrix $\tilde{A}^{(s)} = (a_{ij}^{(s)} - \frac{1}{2}\sigma_i^2)$, which is just the original payoff matrix $A$. Theorem 4.11(a) shows that every strict Nash equilibrium of $A$ is stochastically asymptotically stable for $\{X^{(s)}(t)\}$. This improves a recent result of Khasminskii and Potsepun (2006); see Section 5.

**2. Time averages.** Throughout we assume that $\{X(t)\} = \{X(t,\omega)\}$ is a strong solution to (1.4) with initial value $\mathbf{x}_0 \in \operatorname{int}\Delta$. Then $\{X(t)\}$ is a Markov process and it is not difficult to see that, a.s., $X(t)$ stays in $\operatorname{int}\Delta$ for all $t \geq 0$. Note that if some strategies were not present initially, they would never occur, so that our results would apply to the stochastic dynamics in the interior of the corresponding face of $\Delta$.

We first prove a result that relates time averages of $X(t)$ and Nash equilibria of $\tilde{A}$. The connection with space averages will be considered in Section 3. Recall that $\mathbf{p} \in \Delta$ is a Nash equilibrium of $\tilde{A}$ if

$$\mathbf{p}^T \tilde{A} \mathbf{p} \geq \mathbf{q}^T \tilde{A} \mathbf{p} \qquad \text{for all } \mathbf{q} \in \Delta.$$

The equilibrium is said to be strict if the inequality is strict for all $\mathbf{q} \neq \mathbf{p}$.

THEOREM 2.1.   (a) *For almost every $\omega$ for which*

$$\frac{1}{T}\int_0^T X(t,\omega)\,dt$$

*converges as $T \to \infty$, the limit is a Nash equilibrium of $\tilde{A}$.*

(b) *The following implication holds for almost all $\omega$. If the limit*

$$\mathbf{y} = \mathbf{y}(\omega) = \lim_{k\to\infty} \frac{1}{T_k(\omega)} \int_0^{T_k(\omega)} X(t,\omega)\,dt$$

*exists for a sequence of times $T_k(\omega)$ with*

(2.1)    $T_k(\omega) \nearrow \infty$   *and*   $\dfrac{\log X_i(T_k(\omega),\omega)}{T_k(\omega)} \to 0, \qquad i = 1,\ldots,n,$

*then $\mathbf{y}$ is a Nash equilibrium of $\tilde{A}$ with $\{\tilde{A}\mathbf{y}\}_1 = \cdots = \{\tilde{A}\mathbf{y}\}_n$.*

PROOF. (a) Note first that from (1.3) and Itô's formula,

(2.2) $$d\log Z_i(t) = [\{AX(t)\}_i - \tfrac{1}{2}\sigma_i^2]\,dt + \sigma_i\,dW_i(t),$$

so that

(2.3) $$\begin{aligned}d\log X_i(t) &- d\log X_j(t) \\ &= d\log Z_i(t) - d\log Z_j(t) \\ &= [\{\tilde{A}X(t)\}_i - \{\tilde{A}X(t)\}_j]\,dt + \sigma_i\,dW_i(t) - \sigma_j\,dW_j(t).\end{aligned}$$

The event $\Omega_0 = \{\lim_{t\to\infty} W_j(t)/t \to 0 \text{ for } j=1,\ldots,n\}$ has probability 1. Let $\omega_0 \in \Omega_0$ be such that $\frac{1}{T}\int_0^T X(t,\omega_0)\,dt$ converges to $\mathbf{y} = \mathbf{y}(\omega_0)$, say. Clearly, $\mathbf{y} \in \Delta$. If $y_i > 0$, then there exist times $T_k = T_k(\omega_0)$ such that $0 < T_1 < T_2 < \cdots$, $T_k \to \infty$ and $X_i(T_k,\omega_0) > \frac{1}{2}y_i$ for all $k$. Thus $[\log X_i(T_k,\omega_0)]/T_k \to 0$. As $\omega_0 \in \Omega_0$, it now follows by (2.3) that for every $j$,

$$\begin{aligned}\{\tilde{A}\mathbf{y}\}_i - \{\tilde{A}\mathbf{y}\}_j &= \lim_{k\to\infty} \frac{1}{T_k}\int_0^{T_k} \{\tilde{A}X(t,\omega_0)\}_i - \{\tilde{A}X(t,\omega_0)\}_j\,dt \\ &= \lim_{k\to\infty} \frac{\log X_i(T_k,\omega_0) - \log X_j(T_k,\omega_0)}{T_k} \geq 0.\end{aligned}$$

This shows that $\mathbf{y}$ is a Nash equilibrium of $\tilde{A}$.

(b) The proof is similar to that of (a). □

For every fixed $\omega$ consider the set of limit points $(p_{ij}) = (p_{ij}(\omega))$ of the time averages of the products $X_i(t,\omega)X_j(t,\omega)$ counting the encounters of plays of strategies $i$ against $j$, that is, for some sequence of times $T_k(\omega)$,

(2.4) $$T_k(\omega) \to \infty, \qquad p_{ij} = \lim_{k\to\infty} \frac{1}{T_k(\omega)} \int_0^{T_k(\omega)} X_i(t,\omega)X_j(t,\omega)\,dt \qquad \forall i,j.$$

Let $p_i = \sum_{j=1}^n p_{ij} = \lim_{k\to\infty} \frac{1}{T_k(\omega)} \int_0^{T_k(\omega)} X_i(t,\omega)\,dt$ be its marginals.

THEOREM 2.2. *With probability 1, each limit point $(p_{ij})$ satisfies*

(2.5) $$(\tilde{A}\mathbf{p})_l \leq \sum_{i,j} \tilde{a}_{ij}p_{ij} + \frac{1}{2}\sum_{j=1}^n \sigma_j^2(p_j - p_{jj})$$

*for all $l$, and equality holds for at least one $l$.*

In the deterministic case $\sigma_i = 0$ for all $i$ this result was obtained in Hofbauer (2005). The corresponding inequalities state that the correlated distribution $(p_{ij})$ satisfies the *exact marginal best response property* of Fudenberg and Levine (1995), or is an element of the *Hannan set* as defined by Hart (2005), or a *coarse correlated equilibrium* in the terminology of Young (2004).

PROOF OF THEOREM 2.2. From (1.3) we get for $S = Z_1 + \cdots + Z_n$

$$dS = S\left[X^T A X\, dt + \sum_j \sigma_j X_j\, dW_j\right].$$

Together with $d\log X_l = d\log Z_l - d\log S$ and (2.2) we obtain

(2.6)
$$d\log X_l = \left[(\tilde{A}X)_l - X^T \tilde{A} X - \frac{1}{2}\sum_{j=1}^n \sigma_j^2 X_j(1-X_j)\right] dt + \sigma_l\, dW_l$$
$$- \sum_j \sigma_j X_j\, dW_j.$$

Integrating (2.6) gives

(2.7)
$$\frac{\log X_l(T) - \log X_l(0)}{T} = \frac{1}{T}\int_0^T (\tilde{A}X(t))_l\, dt - \sum_{ij}\tilde{a}_{ij}\frac{1}{T}\int_0^T X_i(t)X_j(t)\, dt$$
$$- \frac{1}{2}\sum_{j=1}^n \sigma_j^2 \frac{1}{T}\int_0^T X_j(t)(1-X_j(t))\, dt + r_l(T)$$

with

$$r_l(T) = \sigma_l \frac{W_l(T)}{T} - \sum_j \sigma_j \frac{1}{T}\int_0^T X_j(t)\, dW_j(t).$$

Since $X_j(t)$ is nonanticipative and $|X_j(t)| \leq 1$, it follows from Friedman [(1975), Corollary 4.6, page 77] that

$$\lim_{T\to\infty}\frac{1}{T}\int_0^T X_j(t)\, dW_j(t) = 0 \qquad \text{a.s.}$$

Consequently,

(2.8)
$$\lim_{T\to\infty} r_l(T) = 0 \qquad \text{a.s.}$$

Since $\log X_l(T) \leq 0$, every limit point of the left-hand side of (2.7) is nonpositive and hence the inequalities (2.5) follow. Now looking at one limit point $(p_{ij})$ arising for the sequence $T_k \to \infty$, let (after possibly refining the time sequence) $\tilde{\mathbf{x}} = \lim_{T_k\to\infty} X(T_k)$ and choose $l$ s.t. $\tilde{x}_l > 0$. Then the left-hand side in (2.7) goes to 0 and there is equality in (2.5). □

REMARK 2.3. For almost every $\omega$ for which the sequence of times $T_k(\omega)$ satisfies (2.1), the corresponding limit matrix $(p_{ij}) = (p_{ij}(\omega))$ satisfies equality in (2.5) for all $l$.

**3. Recurrence.** If the $\omega$-limit of the orbit of $\boldsymbol{\xi}(t)$ given by the deterministic equation (1.2) is contained in $\operatorname{int}\Delta$ and $A$ has a unique interior equilibrium, then the time averages $\frac{1}{T}\int_0^T \boldsymbol{\xi}(t)\,dt$ converge to that equilibrium; see Schuster et al. (1981) or Hofbauer and Sigmund (1998), Theorem 7.6.4. The following theorem contains a stochastic counterpart and also establishes a connection with space averages.

We write $P_{\mathbf{x}_0}$ to indicate the probability computed under the condition $X(0)=\mathbf{x}_0$ and $E_{\mathbf{x}_0}$ to denote expectation with respect to $P_{\mathbf{x}_0}$. For $\mathbf{x}\in\Delta$ and $\varepsilon>0$ let $B_\varepsilon(\mathbf{x})=\{\mathbf{y}\in\Delta:\|\mathbf{y}-\mathbf{x}\|<\varepsilon\}$, where $\|\cdot\|$ denotes the Euclidean norm. The process $\{X(t)\}$ is called recurrent if for every $\mathbf{x}\in\operatorname{int}\Delta$ and every $\varepsilon>0$,

$$P_{\mathbf{x}_0}\{X(T_k)\in B_\varepsilon(\mathbf{x}) \text{ for a sequence of finite random times } T_k$$
$$\text{increasing to infinity}\}=1.$$

See Bhattacharya (1978) for various characterizations of recurrence of diffusions. If $\{X(t)\}$ is recurrent, its transition probability function has a unique (up to positive multiples) $\sigma$-finite invariant measure $\nu$ on $\operatorname{int}\Delta$ [Khas'minskii (1960)]. If $\nu(\operatorname{int}\Delta)<\infty$, $\{X(t)\}$ is said to be positive recurrent. In this case, we will normalize $\nu$ and consider the invariant distribution $\pi=[\nu(\operatorname{int}\Delta)]^{-1}\nu$. If $\{X(t)\}$ is recurrent and $\nu(\operatorname{int}\Delta)=\infty$, $\{X(t)\}$ is said to be null recurrent.

THEOREM 3.1 (Averaging principle). (a) *If* $\{X(t)\}$ *is positive recurrent with invariant distribution* $\pi$ *on* $\operatorname{int}\Delta$, *then* $\hat{\mathbf{y}}:=\int\mathbf{y}\,d\pi(\mathbf{y})$ *is the unique interior Nash equilibrium of* $\tilde{A}$ *and*

$$\lim_{T\to\infty}\frac{1}{T}\int_0^T X(t)\,dt = \hat{\mathbf{y}} \qquad a.s.$$

(b) *Let* $\mathcal{N}=\{\mathbf{y}\in\Delta:\{\tilde{A}\mathbf{y}\}_1=\cdots=\{\tilde{A}\mathbf{y}\}_n\}$. *If* $\{X(t)\}$ *is null recurrent, then* $\mathcal{N}\neq\varnothing$. *Moreover, if* $(T_k)_{k=1}^\infty$ *is a sequence of random times such that*

$$(3.1)\quad T_k\nearrow\infty \quad a.s. \quad and \quad \frac{\log X_i(T_k)}{T_k}\to 0 \qquad a.s.,\ i=1,\ldots,n,$$

*then, a.s., every accumulation point of the sequence*

$$\frac{1}{T_k}\int_0^{T_k} X(t)\,dt, \qquad k=1,2,\ldots,$$

*belongs to the set* $\mathcal{N}$.

PROOF. (a) If $\{X(t)\}$ is positive recurrent with invariant distribution $\pi$ on $\operatorname{int}\Delta$, then $\hat{\mathbf{y}}:=\int\mathbf{y}\,d\pi(\mathbf{y})\in\operatorname{int}\Delta$. Moreover, by the ergodic theorem [Bhattacharya and Waymire (1990), page 623], $\lim_{T\to\infty} \frac{1}{T}\int_0^T X(t)\,dt = \hat{\mathbf{y}}$ a.s. It follows from Theorem 2.1(a) that $\hat{\mathbf{y}}$ is an equilibrium of $\tilde{A}$.

It remains to show that $\hat{\mathbf{y}}$ is the only interior equilibrium of $\tilde{A}$. If not, then there is a line of equilibria which implies the existence of a nonzero vector $\mathbf{c}\in\mathbb{R}^n$ such that $\sum c_i = 0$ and $\mathbf{c}^T\tilde{A}\mathbf{x} = 0$ holds for all $\mathbf{x}\in\Delta$. This together with (2.6) implies

$$(3.2) \qquad \sum_{i=1}^n c_i d\log X_i(t) = \sum_{i=1}^n c_i\sigma_i\,dW_i(t).$$

Integrating gives

$$\sum_{i=1}^n c_i \log \frac{X_i(t)}{X_i(0)} = \sum_{i=1}^n c_i\sigma_i W_i(t).$$

The right-hand side is a one-dimensional Brownian motion and hence null recurrent. This contradicts the assumption that the process $\{X(t)\}$ is positive recurrent in $\operatorname{int}\Delta$.

(b) The assertion is an immediate consequence of Theorem 2.1(b). Note that null recurrence implies that there exists a sequence $(T_k)$ satisfying (3.1). □

REMARK 3.2. (a) Generically, the set $\mathcal{N}$ defined in Theorem 3.1(b) contains at most one element. Thus, if $\mathcal{N} = \{\mathbf{y}\}$, then

$$\lim_{k\to\infty} \frac{1}{T_k}\int_0^{T_k} X(t)\,dt = \mathbf{y} \qquad \text{a.s.}$$

for every sequence of random times $(T_k)$ with (3.1).

(b) Every interior Nash equilibrium of $\tilde{A}$ belongs to $\mathcal{N}$. Thus if $\{X(t)\}$ is positive recurrent, then $\hat{\mathbf{y}}\in\mathcal{N}$, and as $\hat{\mathbf{y}}$ is the only interior Nash equilibrium, $\mathcal{N} = \{\hat{\mathbf{y}}\}$.

THEOREM 3.3. *If $\{X(t)\}$ is positive recurrent with invariant distribution $\pi$ on $\operatorname{int}\Delta$, then the time averages in (2.4) converge a.s. to $p_{ij} := \int y_i y_j\,d\pi(\mathbf{y})$. Therefore, a.s., there is a unique limit matrix, $(p_{ij})$, all its elements are positive, and there holds equality in (2.5) for all $l$.*

PROOF. This follows again from the ergodic theorem. For the last statement, choose a time sequence $T_k\to\infty$ with (3.1) and apply Remark 2.3. □

For $n = 2$ we can compute this limit matrix explicitly: Let $A = \begin{pmatrix} 0 & 1 \\ 1 & 0 \end{pmatrix}$ and $\sigma_1 = \sigma_2 = \sigma > 0$. Then by Proposition 1 of Fudenberg and Harris (1992) or Corollary 3.10, $\{X(t)\}$ is positive recurrent, so that, by Theorem 3.1(a), $p_1 = p_2 = \frac{1}{2}$. Hence $(\tilde{A}\mathbf{p})_1 = \frac{1}{2} - \frac{\sigma^2}{2}$. Since there holds equality in (2.5) and since $\sum_j p_{ij} = p_i$ and $p_{12} = p_{21}$, we obtain

$$(\tilde{A}\mathbf{p})_1 = (p_{12} + p_{21}) - \frac{\sigma^2}{2}(p_{11} + p_{22}) = 2p_{12} - \frac{\sigma^2}{2}(1 - 2p_{12}).$$

It follows that $p_{12} = p_{21} = \frac{1}{4+2\sigma^2} < \frac{1}{4}$ and $p_{11} = p_{22} = \frac{1}{2} - p_{12} > \frac{1}{4}$. For large $\sigma$ the matrix $(p_{ij})$ approaches $\begin{pmatrix} 1/2 & 0 \\ 0 & 1/2 \end{pmatrix}$, which means that the process spends most of the time near the pure strategy states.

REMARK 3.4. If the process $\{X(t)\}$ is null recurrent, then for almost all $\omega$ there is at least one accumulation point $(p_{ij}) = (p_{ij}(\omega))$ for which equality holds for all $l$ in (2.5).

COROLLARY 3.5. *Suppose $\tilde{A} = -\tilde{A}^T$ is a zero-sum game. Then the process $\{X(t)\}$ is not positive recurrent.*

PROOF. If it were positive recurrent, then, by Theorem 3.1, $\tilde{A}\mathbf{p} = 0$ for some $\mathbf{p} \in \operatorname{int} \Delta$, and Theorem 3.3 implies equality for all $l$ in (2.5). Therefore $\sum_{j=1}^n \sigma_j^2(p_j - p_{jj}) = 0$. Since $\sigma_j > 0$ we get $p_j = p_{jj}$ which contradicts $p_{ij} > 0$ for all $i, j$ in Theorem 3.3. □

We next investigate invariant densities and recurrence properties of $\{X(t)\}$ for payoff matrices that satisfy for some $\gamma \in \mathbb{R}$ the condition

$$(3.3) \qquad a_{ij} + a_{ji} - a_{ii} - a_{jj} = \frac{\gamma}{2}(\sigma_i^2 + \sigma_j^2) \qquad \text{for all } i \neq j.$$

If $n = 2$, (3.3) is always satisfied for a certain $\gamma$. For general $n$, the payoff matrices satisfying (3.3) can be characterized as those obtained by subtracting diagonal elements $\gamma\sigma_j^2/2$ from the payoff matrix of a zero-sum game, and rescaling [i.e., adding multiples of $\mathbf{1} = (1, \ldots, 1)^T$ to each column] or adding multiples of $\mathbf{1}^T$ to each row. Note that condition (3.3) is equivalent to the condition obtained by replacing every $a_{kl}$ in (3.3) by $\tilde{a}_{kl}$.

THEOREM 3.6. *Let $\boldsymbol{\alpha} \in \mathbb{R}^n$ and set $\gamma = \alpha_1 + \cdots + \alpha_n$. The transition probability function of $\{X(t)\}$ has an invariant density (with respect to Lebesgue measure on $\operatorname{int} \Delta$) given by*

$$\prod_{i=1}^n x_i^{\alpha_i - 1}, \qquad (x_1, \ldots, x_n) \in \operatorname{int} \Delta,$$

*if and only if one of the following two conditions* (a) *and* (b) *holds.*

(a) *$A$ satisfies (3.3) with $\gamma \neq -1$ and $\{\tilde{A}\boldsymbol{\alpha}\}_1 = \cdots = \{\tilde{A}\boldsymbol{\alpha}\}_n$ and if $\gamma = 0$, then $\{\tilde{A}\boldsymbol{\alpha}\}_1 = \sum_{i=1}^n \tilde{a}_{ii}\alpha_i - \frac{1}{2}\sum_{i=1}^n \sigma_i^2 \alpha_i^2$.*

(b) $\gamma = -1$ and for $k = 1, \ldots, n$,

$$\frac{\sigma_k^2}{2} + \frac{1}{2}\sum_{i=1}^{n}\sigma_i^2\alpha_i^2 + \alpha_k\sigma_k^2 - \tilde{a}_{kk} - \sum_{i=1}^{n}\alpha_i\tilde{a}_{ik} = 0.$$

PROOF. For simplicity, we work with the $\mathbb{R}^{n-1}$-valued process $Y(t) = \Psi(X(t))$, where $\Psi(x_1, \ldots, x_n) = (\log(x_1/x_n), \ldots, \log(x_{n-1}/x_n))^T$. Let $\mathbf{e}_k$ denote the $k$th column of the $n \times n$ identity matrix. We have from (2.3)

$$dY_j(t) = \{(\mathbf{e}_j - \mathbf{e}_n)^T \tilde{A}\Psi^{-1}(Y(t))\}\,dt + \sigma_j\,dW_j(t) - \sigma_n\,dW_n(t),$$

$j = 1, \ldots, n-1$. Let $\tilde{L}$ denote the formal adjoint of the differential operator corresponding to $\{Y(t)\}$, that is,

$$\tilde{L}f(\mathbf{y}) = \frac{1}{2}\sum_{i,j=1}^{n-1}\frac{\partial^2}{\partial y_i\,\partial y_j}[(\sigma_n^2 + \delta_{ij}\sigma_i^2)f(\mathbf{y})] - \sum_{i=1}^{n-1}\frac{\partial}{\partial y_i}[(\mathbf{e}_i - \mathbf{e}_n)^T\tilde{A}\Psi^{-1}(y)f(\mathbf{y})]$$

for $f \in C^2(\mathbb{R}^{n-1})$. A change of variables shows that $\prod_{i=1}^{n} x_i^{\alpha_i - 1}$ is an invariant density for the transition probability function of $\{X(t)\}$ if and only if

$$g(\mathbf{y}) = \frac{\exp(\sum_{i=1}^{n-1}\alpha_i y_i)}{(1 + e^{y_1} + \cdots + e^{y_{n-1}})^{\sum_{i=1}^{n}\alpha_i}}, \qquad \mathbf{y} \in \mathbb{R}^{n-1},$$

is an invariant density for the transition probability function of $\{Y(t)\}$.

Set

$$\phi(\mathbf{y}) = 1 + e^{y_1} + \cdots + e^{y_{n-1}}.$$

Then for all $i, j = 1, \ldots, n$, $i \neq j$,

$$\frac{\partial g(\mathbf{y})}{\partial y_i} = \frac{g(\mathbf{y})}{\phi(\mathbf{y})}[\alpha_i\phi(\mathbf{y}) - \gamma e^{y_i}],$$

$$\frac{\partial^2 g(\mathbf{y})}{\partial y_i^2} = \frac{g(\mathbf{y})}{\phi^2(\mathbf{y})}\{[\alpha_i\phi(\mathbf{y}) - \gamma e^{y_i}]^2 + \gamma e^{y_i}[e^{y_i} - \phi(\mathbf{y})]\},$$

$$\frac{\partial^2 g(\mathbf{y})}{\partial y_i\,\partial y_j} = \frac{g(\mathbf{y})}{\phi^2(\mathbf{y})}\{[\alpha_i\phi(\mathbf{y}) - \gamma e^{y_i}][\alpha_j\phi(\mathbf{y}) - \gamma e^{y_j}] + \gamma e^{y_i + y_j}\}.$$

Hence

$$(3.4)\quad \tilde{L}g(\mathbf{y}) = \frac{g(\mathbf{y})}{\phi^2(\mathbf{y})}\left\{\sum_{1 \leq k < l \leq n-1}\rho_{kl}e^{y_k + y_l} + \sum_{k=1}^{n-1}\rho_{kn}e^{y_k} + \rho_n + \sum_{k=1}^{n-1}\rho_k e^{2y_k}\right\},$$

where

$$\rho_k = \frac{1}{2}\left\{\sum_{i=1}^{n}\sigma_i^2\alpha_i^2 + \gamma\sigma_k^2(\gamma - 2\alpha_k)\right\} + \sum_{i=1}^{n}\alpha_i(\tilde{a}_{kk} - \tilde{a}_{ik}),$$

$$\rho_{kl} = \sum_{i=1}^{n} \sigma_i^2 \alpha_i^2 - \gamma \left[ \sigma_k^2 \left( \alpha_k + \frac{1}{2} \right) + \sigma_l^2 \left( \alpha_l + \frac{1}{2} \right) \right]$$

$$- \tilde{a}_{kk} - \tilde{a}_{ll} + (\gamma + 1)(\tilde{a}_{kl} + \tilde{a}_{lk}) - \sum_{i=1}^{n} \alpha_i (\tilde{a}_{ik} + \tilde{a}_{il}).$$

For $k \neq l$,

$$(3.5) \quad \rho_{kl} - \rho_k - \rho_l = (\gamma + 1) \left[ a_{kl} + a_{lk} - a_{kk} - a_{ll} - \frac{\gamma}{2} (\sigma_k^2 + \sigma_l^2) \right].$$

Furthermore,

$$(3.6) \quad \rho_1 = \cdots = \rho_n = 0 \iff B\boldsymbol{\alpha} = \left( \frac{1}{2} \sum_{i=1}^{n} \sigma_i^2 \alpha_i^2 \right) \mathbf{1},$$

where

$$B = \tilde{A}^T - \tilde{\mathbf{a}} \mathbf{1}^T - \frac{\gamma}{2} \mathbf{s} \mathbf{1}^T + \gamma \operatorname{diag}(\sigma_1^2, \ldots, \sigma_n^2),$$

$$\tilde{\mathbf{a}} = (\tilde{a}_{11}, \ldots, \tilde{a}_{nn})^T, \qquad \mathbf{s} = (\sigma_1^2, \ldots, \sigma_n^2)^T.$$

Condition (3.3) is equivalent to

$$(3.7) \quad B = -\tilde{A} + \mathbf{1} \tilde{\mathbf{a}}^T + \frac{\gamma}{2} \mathbf{1} \mathbf{s}^T.$$

Now suppose $\prod_{i=1}^{n} x_i^{\alpha_i - 1}$ is an invariant density for the transition probability function of $\{X(t)\}$. Then $\widetilde{L}g = 0$; see Theorem 8.4 in Pinsky (1995), page 181. Hence, in view of (3.4), $\rho_k = 0$ for all $k$ and $\rho_{kl} = 0$ for all $k \neq l$. Suppose first that $\gamma \neq -1$. It then follows by (3.5) that $A$ satisfies (3.3). Moreover, by (3.6) and (3.7),

$$\left( \frac{1}{2} \sum_{i=1}^{n} \sigma_i^2 \alpha_i^2 \right) \mathbf{1} = B\boldsymbol{\alpha} = \left( \tilde{\mathbf{a}}^T \boldsymbol{\alpha} + \frac{\gamma}{2} \mathbf{s}^T \boldsymbol{\alpha} \right) \mathbf{1} - \tilde{A} \boldsymbol{\alpha}.$$

This shows that $\tilde{A}\boldsymbol{\alpha}$ is proportional to $\mathbf{1}$, and if $\gamma = 0$, then $\{\tilde{A}\boldsymbol{\alpha}\}_1 = \sum_{i=1}^{n} \tilde{a}_{ii} \alpha_i - \frac{1}{2} \sum_{i=1}^{n} \sigma_i^2 \alpha_i^2$. Thus condition (a) is satisfied. If $\gamma = -1$, then condition (b) holds because $\rho_k = 0$ for all $k$.

Conversely, suppose condition (a) holds. Then, by (3.7),

$$(3.8) \quad B\boldsymbol{\alpha} = \left( \tilde{\mathbf{a}}^T \boldsymbol{\alpha} + \frac{\gamma}{2} \mathbf{s}^T \boldsymbol{\alpha} \right) \mathbf{1} - \tilde{A}\boldsymbol{\alpha} = \mu \mathbf{1}$$

for some $\mu \in \mathbb{R}$. Thus $\boldsymbol{\alpha}^T B \boldsymbol{\alpha} = \gamma \mu$. By (3.7), $B^T + B = \gamma \operatorname{diag}(\sigma_1^2, \ldots, \sigma_n^2)$, and so $2\gamma\mu = \boldsymbol{\alpha}^T (B^T + B) \boldsymbol{\alpha} = \gamma \sum_{i=1}^{n} \alpha_i^2 \sigma_i^2$. If $\gamma \neq 0$, it follows that $\mu = \frac{1}{2} \sum_{i=1}^{n} \alpha_i^2 \sigma_i^2$. If $\gamma = 0$, then, according to condition (a), $\tilde{A}\boldsymbol{\alpha} = [(\sum_{i=1}^{n} \alpha_i (\tilde{a}_{ii} - \frac{1}{2} \sigma_i^2 \alpha_i)] \mathbf{1}$, and it is obvious from (3.8) that again $\mu = \frac{1}{2} \sum_{i=1}^{n} \alpha_i^2 \sigma_i^2$. It therefore follows from (3.8) and (3.6) that $\rho_k = 0$ for all $k$. Equations (3.3) and (3.5) now imply that $\rho_{kl} = 0$ for all $k \neq l$. Consequently, by (3.4), $\widetilde{L}g = 0$. If condition (b) holds, then $\rho_k = 0$ for all $k$, and by (3.5), $\rho_{kl} = 0$ for all $k \neq l$. Thus, again, $\widetilde{L}g = 0$. It now follows by Theorem 8.5 in Pinsky [(1995), page 182] that $g$ is an invariant density, provided the diffusion corresponding to the solution to the generalized martingale problem for $\widetilde{L}^g$ on $\mathbb{R}^{n-1}$ does not explode, where $\widetilde{L}^g f = \frac{1}{g}\widetilde{L}(gf)$. We have

$$\widetilde{L}^g f(\mathbf{y}) = \widetilde{L}_0 f(\mathbf{y}) + \sum_{j=1}^{n-1}\left\{\sum_{i=1}^{n-1} \frac{(\sigma_n^2 + \delta_{ij}\sigma_i^2)}{g(\mathbf{y})}\frac{\partial g(\mathbf{y})}{y_i}\right\}\frac{\partial f(\mathbf{y})}{\partial y_j},$$

where

$$\widetilde{L}_0 f(\mathbf{y}) = \frac{1}{2}\sum_{i,j=1}^{n-1}\frac{\partial^2}{\partial y_i\,\partial y_j}[(\sigma_n^2 + \delta_{ij}\sigma_i^2)f(\mathbf{y})] - \sum_{i=1}^{n-1}(\mathbf{e}_i - \mathbf{e}_n)^T \tilde{A}\Psi^{-1}(y)\frac{\partial f(\mathbf{y})}{\partial y_i}$$

and

$$\frac{1}{g(\mathbf{y})}\frac{\partial g(\mathbf{y})}{y_i} = \alpha_i - \frac{\gamma e^{y_i}}{\phi(\mathbf{y})}.$$

Thus the coefficients of $\widetilde{L}^g$ are Lipschitz continuous and bounded so that the martingale problem for $\widetilde{L}^g$ has a unique solution and explosions do not occur. $\square$

REMARK 3.7. Suppose $A$ satisfies (3.3). Then $A$ is conditionally negative [positive] definite if $\gamma > [<]0$. Indeed, if $y_1 + \cdots + y_n = 0$, then, by (3.3),

$$2\mathbf{y}^T A\mathbf{y} = \sum_{i,j:\,i\neq j} y_i(a_{ij} + a_{ji} - a_{ii} - a_{jj})y_j = -\gamma \sum_i \sigma_i^2 y_i^2.$$

In particular, if (3.3) holds with $\gamma \neq 0$, then there can be at most one $\boldsymbol{\alpha}$ with $\alpha_1 + \cdots + \alpha_n = \gamma$ that satisfies the equations $\{\tilde{A}\boldsymbol{\alpha}\}_1 = \cdots = \{\tilde{A}\boldsymbol{\alpha}\}_n$ in condition (a) of Theorem 3.6.

The density $\prod_{i=1}^n x_i^{\alpha_i - 1}$ in Theorem 3.6 is integrable if and only if $\alpha_i > 0$ for all $i$. A random variable $U$ with $P\{U \in \text{int}\,\Delta\} = 1$ which has a density that is proportional to $\prod_{i=1}^n x_i^{\alpha_i - 1}$ is said to have a Dirichlet distribution with parameter $\boldsymbol{\alpha} = (\alpha_1, \ldots, \alpha_n)^T$, provided $\alpha_i > 0$ for all $i$. If $U$ has such a distribution and $\gamma = \alpha_1 + \cdots + \alpha_n$, then

(3.9) $$\mathsf{E}U = \gamma^{-1}\boldsymbol{\alpha}, \qquad \text{Var}[U_i] = \frac{\alpha_i(\gamma - \alpha_i)}{\gamma^2(\gamma + 1)}.$$

COROLLARY 3.8. *The process $\{X(t)\}$ is positive recurrent and its transition probability function has an invariant Dirichlet distribution with parameter $\boldsymbol{\alpha}$ if and only if $A$ satisfies (3.3) with some $\gamma > 0$ and $\gamma^{-1}\boldsymbol{\alpha}$ is an interior Nash equilibrium of $\tilde{A}$.*

REMARK 3.9. If $n = 2$, then the condition that $A$ satisfies (3.3) with some $\gamma > 0$ and $\tilde{A}$ has an interior Nash equilibrium is equivalent to the inequalities $\tilde{a}_{12} > \tilde{a}_{22}$ and $\tilde{a}_{21} > \tilde{a}_{11}$. For this case, the stable coexistence case, Fudenberg and Harris (1992) already calculated the ergodic distribution of $\{X(t)\}$ by means of the normalized speed measure.

COROLLARY 3.10. *Suppose $A$ satisfies (3.3) for some $\gamma \in \mathbb{R}$. Then $\{X(t)\}$ is positive recurrent if and only if $\tilde{A}$ has an interior Nash equilibrium and $\gamma > 0$.*

PROOF. To prove the sufficiency of the condition suppose that $\mathbf{x}$ is an interior Nash equilibrium of $\tilde{A}$ and $\gamma > 0$. Set $\boldsymbol{\alpha} = \gamma \mathbf{x}$. Then $\gamma^{-1}\boldsymbol{\alpha}$ is an interior Nash equilibrium and so, by Corollary 3.8, $\{X(t)\}$ is positive recurrent.

To prove necessity, suppose $\{X(t)\}$ is positive recurrent. Then, by Theorem 3.1(a), $\tilde{A}$ has an interior Nash equilibrium $\mathbf{x}$, say. Let $\boldsymbol{\alpha} = \gamma \mathbf{x}$. If $\gamma \neq -1$, then Theorem 3.6 shows that $\prod_{i=1}^{n} x_i^{\alpha_i - 1}$ is an invariant density, which must be integrable, so that $\gamma > 0$. It remains to rule out that $\gamma = -1$. If $\gamma = -1$, $A$ would be conditionally positive definite by Remark 3.7, so that Theorem 4.5 would yield that $\{X(t)\}$ is transient. □

REMARK 3.11. Corollaries 3.8 and 3.10 show that for games with (3.3) only Dirichlet distributions can be invariant distributions, and if $\{X(t)\}$ is positive recurrent and $\sigma_1 = \cdots = \sigma_n$, then $\{X(t)\}$ remains positive recurrent if $\sigma_1, \ldots, \sigma_n$ are replaced with $\kappa\sigma_1, \ldots, \kappa\sigma_n$ for any $\kappa > 0$. As $\kappa \searrow 0$, $\gamma \nearrow \infty$ and, in view of (3.9), the invariant distribution converges weakly to the unit mass concentrated at the Nash equilibrium of $A$. As $\kappa \nearrow \infty$, the invariant distribution converges vaguely to the zero measure on $\operatorname{int}\Delta$.

COROLLARY 3.12. *The density $\prod_{i=1}^{n} x_i^{-1}$, $\mathbf{x} \in \operatorname{int}\Delta$, is invariant for the transition probability function of $\{X(t)\}$ if and only if*

$$(3.10) \qquad a_{ij} + a_{ji} - a_{ii} - a_{jj} = 0 \quad \text{for all } i,j.$$

*Moreover, if (3.10) holds and there exists $\boldsymbol{\beta} \in \mathbb{R}^n$ such that*

$$(3.11) \quad \beta_1 + \cdots + \beta_n = 0, \qquad \{\tilde{A}\boldsymbol{\beta}\}_1 = \cdots = \{\tilde{A}\boldsymbol{\beta}\}_n \quad \text{and} \quad \boldsymbol{\beta}^T \tilde{A} \neq \mathbf{0},$$

*then $\prod_{i=1}^{n} x_i^{\alpha_i - 1}$ is another invariant density, where $\alpha_i = c\beta_i$ and $c = 2\{\boldsymbol{\beta}^T \tilde{A}\}_1 / (\sum_{i=1}^{n} \sigma_i^2 \beta_i^2) \neq 0$. In this case, $\{X(t)\}$ is transient.*

PROOF. The claimed equivalence is an immediate consequence of Theorem 3.6. If (3.10) and (3.11) hold, then

$$\sum_{i=1}^{n} \tilde{a}_{ii}\alpha_i - \frac{1}{2}\sum_{i=1}^{n}\sigma_i^2\alpha_i^2 = c\sum_{i=1}^{n}(\tilde{a}_{ii} - \tilde{a}_{i1})\beta_i = c\sum_{i=1}^{n}(a_{1i} - a_{11})\beta_i = \{\tilde{A}\boldsymbol{\alpha}\}_1,$$

so that $\prod_{i=1}^{n} x_i^{\alpha_i-1}$ is an invariant density by Theorem 3.6. Furthermore, (3.10) and (3.11) imply that $\{\boldsymbol{\beta}^T \tilde{A}\}_1 = \cdots = \{\boldsymbol{\beta}^T \tilde{A}\}_n$, and so, as $\boldsymbol{\beta}^T \tilde{A} \neq \mathbf{0}$, we have $c \neq 0$. Thus the invariant densities $\prod_{i=1}^{n} x_i^{-1}$ and $\prod_{i=1}^{n} x_i^{\alpha_i-1}$ are not proportional. This implies that $\{X(t)\}$ is transient; see Khas'minskii (1960) or Pinsky (1995), pages 148 and 181. □

REMARK 3.13. If $\tilde{A}$ is a zero-sum game, then (3.11) is equivalent to $\{\tilde{A}\boldsymbol{\beta}\}_1 = \cdots = \{\tilde{A}\boldsymbol{\beta}\}_n \neq 0$.

**4. Transience and stability of Nash equilibria.** The deterministic replicator dynamics satisfies an important exclusion principle: if the replicator equation (1.2) has no interior rest point, that is, if there does not exist $\mathbf{p} \in \text{int}\,\Delta$ such that $\{A\mathbf{p}\}_1 = \cdots = \{A\mathbf{p}\}_n$, then every orbit converges to the boundary of $\Delta$; see Hofbauer (1981) or Hofbauer and Sigmund (1998), Theorem 7.6.1. The following result for the stochastic replicator dynamics yields convergence to the boundary under a slightly stronger condition on $\tilde{A}$.

The process $\{X(t)\}$ is said to be transient if $P_{\mathbf{x}_0}\{X(t) \to \text{bd}\,\Delta\} = 1$.

THEOREM 4.1 (Exclusion principle). *If there does not exist $\mathbf{p} \in \Delta$ such that $\{\tilde{A}\mathbf{p}\}_1 = \cdots = \{\tilde{A}\mathbf{p}\}_n$, then $\{X(t)\}$ is transient.*

Theorem 4.1 is an immediate consequence of the averaging principle (Theorem 3.1) and the dichotomy between recurrence and transience; see, for example, Theorem 8.1 in Pinsky (1995), page 74. In the following we present an *alternative proof*.

PROOF OF THEOREM 4.1. If there is no $\mathbf{p} \in \Delta$ such that $\{\tilde{A}\mathbf{p}\}_1 = \cdots = \{\tilde{A}\mathbf{p}\}_n$, then by a simple separation argument [see, e.g., Hofbauer (1981)] there exists a $\mathbf{c} \in \mathbb{R}^n$ such that $c_1 + \cdots + c_n = 0$ and $\mathbf{c}^T \tilde{A}\mathbf{x} > 0$ for all $\mathbf{x} \in \Delta$. By compactness there is $\delta > 0$ such that $\mathbf{c}^T \tilde{A}\mathbf{x} \geq \delta > 0$ for all $\mathbf{x} \in \Delta$. Since $c_1 + \cdots + c_n = 0$, (2.6) implies

$$(4.1) \qquad \sum_{i=1}^{n} c_i d\log X_i(t) = \mathbf{c}^T \tilde{A} X(t)\,dt + \sum_{i=1}^{n} c_i \sigma_i \, dW_i(t).$$

Integrating gives

$$\sum_{i=1}^{n} c_i \log \frac{X_i(t)}{X_i(0)} \geq \delta t + \sum_{i=1}^{n} c_i \sigma_i W_i(t),$$

and hence $\sum_{i=1}^{n} c_i \log X_i(t) \to \infty$ a.s., as $t \to \infty$, which shows transience. □

COROLLARY 4.2. *If strategy $k$ is strictly dominated in the game $\tilde{A}$ [i.e., there is a strategy $\mathbf{q} \in \Delta$ such that $(\tilde{A}\mathbf{x})_k < \mathbf{q}^T \tilde{A}\mathbf{x}$ for all $\mathbf{x} \in \Delta$], then $X_k(t) \to 0$ a.s.*

PROOF. Choose $\mathbf{c} = \mathbf{q} - \mathbf{e}_k$ in the previous proof which shows

$$\sum_{i=1}^n q_i \log X_i(t) - \log X_k(t) \to \infty \qquad \text{a.s.}$$

Since the sum is bounded above, $\log X_k(t) \to -\infty$ a.s., or $X_k(t) \to 0$ a.s.
□

The following example shows that the weaker condition that $\tilde{A}$ has no interior Nash equilibrium (or that strategy $k$ is weakly dominated) is not sufficient for transience.

EXAMPLE 4.3. Suppose $n = 2$ and $\tilde{a}_{11} = \tilde{a}_{21}$, $\tilde{a}_{12} > \tilde{a}_{22}$. Then $\tilde{A}$ does not have an interior Nash equilibrium. Using Feller's criterion [see, e.g., Proposition 5.22 in Karatzas and Shreve (1991), page 345] one may show that $\{X(t)\}$ is recurrent. In view of Theorem 3.1(a), $\{X(t)\}$ must be null recurrent.

REMARK 4.4. If $A$ satisfies (3.10) and there exists $\boldsymbol{\beta}$ satisfying (3.11), then $\beta_1 + \cdots + \beta_n = 0$ and $\{\boldsymbol{\beta}^T \tilde{A}\}_1 = \cdots = \{\boldsymbol{\beta}^T \tilde{A}\}_n \neq 0$. That $\{X(t)\}$ must then be transient can also be shown as in the alternative proof of Theorem 4.1.

For the second sufficient transience criterion recall that an $n \times n$ matrix $B$ is said to be conditionally positive definite if $\mathbf{y}^T B \mathbf{y} > 0$ for all $\mathbf{y} \neq \mathbf{0}$ with $y_1 + \cdots + y_n = 0$. Note that $\tilde{A}$ is conditionally positive definite if and only if $A$ is conditionally positive definite.

THEOREM 4.5. *Suppose there exists $\mathbf{p} \in \mathbb{R}^n \setminus \{\mathbf{e}_1, \ldots, \mathbf{e}_n\}$ such that $\{\tilde{A}\mathbf{p}\}_1 = \cdots = \{\tilde{A}\mathbf{p}\}_n$ and $p_1 + \cdots + p_n = 1$. Suppose further that $\tilde{A}$ (or $A$) is conditionally positive definite. Then $\min\{X_i(t) : i \in I^+\} \to 0$ a.s., where $I^+ = \{i : p_i > 0\}$. In particular, $\{X(t)\}$ is transient.*

PROOF. Set

$$\phi(\mathbf{x}) = \sum_{j=1}^n p_j \log x_j, \qquad \mathbf{x} \in \operatorname{int} \Delta.$$

By (2.6),

$$\phi(X(t)) = \phi(X(0)) + \int_0^t h(X(s)) \, ds + R(t),$$

where

$$h(\mathbf{x}) = \mathbf{p}^T \tilde{A}\mathbf{x} - \mathbf{x}^T \tilde{A}\mathbf{x} - \frac{1}{2}\sum_{j=1}^n \sigma_j^2 x_j(1-x_j)$$

$$= -(\mathbf{p}-\mathbf{x})^T \tilde{A}(\mathbf{p}-\mathbf{x}) - \frac{1}{2}\sum_{j=1}^n \sigma_j^2 x_j(1-x_j)$$

and

$$R(t) = \sum_{j=1}^n p_j \sigma_j W_j(t) - \sum_{j=1}^n \sigma_j \int_0^t X_j(s)\,dW_j(s).$$

Since $\tilde{A}$ is conditionally positive definite, $(\mathbf{p}-\mathbf{x})^T \tilde{A}(\mathbf{p}-\mathbf{x}) \geq 0$ for all $\mathbf{x} \in \Delta$, with equality if and only if $\mathbf{x} = \mathbf{p}$. Moreover, $\sum_{j=1}^n \sigma_j^2 x_j(1-x_j) \geq 0$ for all $\mathbf{x} \in \Delta$, with equality if and only if $\mathbf{x} \in \{\mathbf{e}_1, \ldots, \mathbf{e}_n\}$. As $\mathbf{p} \in \mathbb{R}^n \setminus \{\mathbf{e}_1, \ldots, \mathbf{e}_n\}$, it follows that there exists some $\varepsilon > 0$ such that

$$h(\mathbf{x}) \leq -\varepsilon \qquad \text{for all } \mathbf{x} \in \operatorname{int} \Delta.$$

In view of (2.8), $R(t)/t \to 0$ a.s. Consequently, a.s.,

$$\limsup_{t\to\infty} \frac{1}{t}\sum_{i\in I^+} p_i \log X_i(t) \leq \limsup_{t\to\infty} \frac{1}{t}\phi(X(t)) \leq -\varepsilon + \limsup_{t\to\infty} \frac{1}{t}R(t) = -\varepsilon.$$

Thus $\min\{X_i(t): i \in I^+\} \to 0$ a.s. $\square$

REMARK 4.6. If $A$ is conditionally positive definite and there exists $\mathbf{p} \in \mathbb{R}^n \setminus \Delta$ such that $\{\tilde{A}\mathbf{p}\}_1 = \cdots = \{\tilde{A}\mathbf{p}\}_n$ and $p_1 + \cdots + p_n = 1$, then there cannot exist $\mathbf{q} \in \Delta$ with $\{\tilde{A}\mathbf{q}\}_1 = \cdots = \{\tilde{A}\mathbf{q}\}_n$, and it follows already by Theorem 4.1 that $X(t)$ converges to $\operatorname{bd}\Delta$ a.s. Theorem 4.5 provides in this case a proper subset of $\operatorname{bd}\Delta$ to which $X(t)$ converges a.s.

REMARK 4.7. Theorem 4.5 complements the following stability result of Imhof (2005a): If $A$ has an interior Nash equilibrium and $A$ is conditionally negative definite, then $\{X(t)\}$ is positive recurrent and spends most of the time in a small neighborhood of the equilibrium, provided $\sigma_1, \ldots, \sigma_n$ are sufficiently small. Note that the transience result holds without restriction on the size of noise.

COROLLARY 4.8. *If $A$ is conditionally positive definite and none of the columns of $\tilde{A}$ is proportional to $\mathbf{1}$, then $\{X(t)\}$ is transient.*

PROOF. If there does not exist $\mathbf{p} \in \Delta$ with $\{\tilde{A}\mathbf{p}\}_1 = \cdots = \{\tilde{A}\mathbf{p}\}_n$, then Theorem 4.1 yields transience. Otherwise, there does exist such a $\mathbf{p} \in \Delta$, but if none of the columns of $\tilde{A}$ is proportional to $\mathbf{1}$, then $\mathbf{p}$ cannot belong to $\{\mathbf{e}_1, \ldots, \mathbf{e}_n\}$. Transience now follows from Theorem 4.5. $\square$

REMARK 4.9. *If one of the columns of $\tilde{A}$ is proportional to $\mathbf{1}$, then $\{X(t)\}$ cannot be positive recurrent.* This follows from Theorem 3.1(a).

EXAMPLE 4.10. In the case $n = 2$ the assumptions of Theorem 4.5 are satisfied if $\tilde{a}_{11} > \tilde{a}_{21}$ and $\tilde{a}_{12} < \tilde{a}_{22}$. This means that the game $\tilde{A}$ is bistable: both pure strategies are strict equilibria.

A point $\mathbf{p} \in \Delta$ is said to be *stochastically asymptotically stable* [Arnold (1974), page 181] if for every open neighborhood $U$ of $\mathbf{p}$ and every $\varepsilon > 0$ there is a neighborhood $V$ of $\mathbf{p}$ such that for every initial state $\mathbf{x} \in V \cap \text{int}\,\Delta$,

$$P_\mathbf{x}\Big\{X(t) \in U \text{ for all } t \geq 0, \lim_{t\to\infty} X(t) = \mathbf{p}\Big\} \geq 1 - \varepsilon.$$

For every $\mathbf{p} \in \Delta \setminus \{\mathbf{e}_1, \ldots, \mathbf{e}_n\}$ and every $\mathbf{x} \in \text{int}\,\Delta$, $P_\mathbf{x}\{X(t) \to \mathbf{p}\} = 0$; see Imhof (2005a). Thus only the vertices of $\Delta$ can be stochastically asymptotically stable.

THEOREM 4.11. (a) *Suppose strategy $k$ is a strict Nash equilibrium in the game $\tilde{A}$, that is,*

$$a_{kk} > a_{jk} + \frac{\sigma_k^2}{2} - \frac{\sigma_j^2}{2} \qquad \text{for all } j \neq k.$$

*Then $\mathbf{e}_k$ is stochastically asymptotically stable.*

(b) *Suppose that for some $\mathbf{x} \in \text{int}\,\Delta$,*

$$P_\mathbf{x}\Big\{\lim_{t\to\infty} X(t) = \mathbf{e}_k\Big\} > 0.$$

*Then strategy $k$ must be a Nash equilibrium of $\tilde{A}$. Moreover, if there exists $i \neq k$ with $\tilde{a}_{ik} = \tilde{a}_{kk}$, then there exists $j \neq k$ such that $\tilde{a}_{ij} < \tilde{a}_{kj}$.*

PROOF. (a) Let $L$ denote the second-order differential operator associated with $\{X(t)\}$, that is,

$$Lf(\mathbf{x}) = \sum_{j=1}^n b_j(\mathbf{x}) \frac{\partial f(\mathbf{x})}{\partial x_j} + \frac{1}{2} \sum_{j,k=1}^n \gamma_{jk}(\mathbf{x}) \frac{\partial^2 f(\mathbf{x})}{\partial x_j \, \partial x_k}, \qquad f \in C^2(\text{int}\,\Delta),$$

where

$$b_j(\mathbf{x}) = x_j(\mathbf{e}_j - \mathbf{x})^T[A - \text{diag}(\sigma_1^2, \ldots, \sigma_n^2)]\mathbf{x},$$

$$\gamma_{jk}(\mathbf{x}) = \sum_{\nu=1}^n c_{j\nu}(\mathbf{x}) c_{k\nu}(\mathbf{x}), \qquad c_{jk}(\mathbf{x}) = \begin{cases} x_j(1-x_j)\sigma_j, & j=k, \\ -x_j x_k \sigma_k, & j \neq k. \end{cases}$$

Consider the Lyapunov function

$$\phi(\mathbf{x}) = \sum_{j \neq k} x_j^r,$$

where $r > 0$ will be chosen later. Clearly, $\phi(\mathbf{x}) \geq 0$ for all $\mathbf{x} \in \Delta$, and $\phi(\mathbf{x}) = 0$ if and only if $\mathbf{x} = \mathbf{e}_k$. We will show that there exists a neighborhood $V$ of $\mathbf{e}_k$ such that for some $c > 0$,

$$L\phi(\mathbf{x}) \leq -c\phi(\mathbf{x}) \qquad \text{for all } \mathbf{x} \in V \cap \text{int}\,\Delta.$$

This yields the assertion by Theorem 4 and Remark 2 in Gichman and Skorochod (1971), pages 314–315, or Theorem 4.1 in Has'minskiĭ (1980), page 167.

For all $\mathbf{x} \in \text{int}\,\Delta$,

$$L\phi(\mathbf{x}) = r\sum_{j \neq k} x_j^r (\mathbf{e}_j - \mathbf{x})^T [A - \text{diag}(\sigma_1^2, \ldots, \sigma_n^2)]\mathbf{x} + \frac{r(r-1)}{2}\sum_{j \neq k} \gamma_{jj}(\mathbf{x}) x_j^{r-2}$$

$$= r\sum_{j \neq k} x_j^r \left\{ (\mathbf{e}_j - \mathbf{x})^T \tilde{A}\mathbf{x} + \frac{r}{2}\left( \sigma_j^2(1 - 2x_j) + \sum_{\nu=1}^n \sigma_\nu^2 x_\nu^2 \right) \right.$$

$$\left. + \frac{1}{2}\sum_{\nu=1}^n \sigma_\nu^2 x_\nu(x_\nu - 1) \right\}$$

$$\leq r\sum_{j \neq k} x_j^r \{(\mathbf{e}_j - \mathbf{x})^T \tilde{A}\mathbf{x} + r\max\{\sigma_1^2, \ldots, \sigma_n^2\}\}.$$

Since $k$ is a strict Nash equilibrium for $\tilde{A}$, there exists $\varepsilon > 0$ such that for all $j \neq k$, $(\mathbf{e}_k - \mathbf{e}_j)^T \tilde{A}\mathbf{e}_k > \varepsilon$. Thus there is a neighborhood $V$ of $\mathbf{e}_k$ such that

$$(\mathbf{x} - \mathbf{e}_j)^T \tilde{A}\mathbf{x} > \frac{\varepsilon}{2} \qquad \text{for all } \mathbf{x} \in V \text{ and } j \neq k.$$

Now if $r > 0$ is so small that $r\max\{\sigma_1^2, \ldots, \sigma_n^2\} < \frac{1}{4}\varepsilon$, then for all $\mathbf{x} \in V \cap \text{int}\,\Delta$

$$L\phi(\mathbf{x}) \leq r\sum_{j \neq k} x_j^r \left( -\frac{\varepsilon}{2} + \frac{\varepsilon}{4} \right) = -\frac{r\varepsilon}{4}\phi(\mathbf{x}).$$

(b) If $\lim_{t \to \infty} X(t, \omega) = \mathbf{e}_k$, then the time averages $\frac{1}{T}\int_0^T X(t, \omega)\,dt$ converge to $\mathbf{e}_k$, too. It follows from Theorem 2.1 that strategy $k$ is a Nash equilibrium. For an indirect proof of the remaining part assume there exists $i \neq k$ such that $\tilde{a}_{ik} = \tilde{a}_{kk}$ and $\tilde{a}_{ij} \geq \tilde{a}_{kj}$ for all $j$. Since two-dimensional Brownian motion is recurrent, there is a sequence of random times $(T_l)_{l=1}^\infty$ such that $T_l \nearrow \infty$ and, for every $l$, $W_i(T_l) \geq 0$ and $W_k(T_l) \leq 0$. It follows that for every $l$,

$$\log\frac{X_i(T_l)}{X_k(T_l)} - \log\frac{X_i(0)}{X_k(0)}$$

$$= \int_0^{T_l} \{\tilde{A}X(t)\}_i - \{\tilde{A}X(t)\}_k\,dt + \sigma_i W_i(T_l) - \sigma_k W_k(T_l) \geq 0,$$

which contradicts the assumption that $P_\mathbf{x}\{\lim_{t \to \infty} X(t) = \mathbf{e}_k\} > 0$. $\square$

REMARK 4.12. Theorem 4.11(a) improves Theorem 4.1 of Imhof (2005a), where the same stability assertion was proved under the stronger condition that

$$a_{kk} > a_{jk} + \sigma_k^2 \qquad \text{for all } j \neq k.$$

Note that in view of Example 4.3 or Theorem 4.11(b), the stability condition in Theorem 4.11(a) cannot be weakened to

$$a_{kk} \geq a_{jk} + \frac{\sigma_k^2}{2} - \frac{\sigma_j^2}{2} \qquad \text{for all } j \neq k.$$

REMARK 4.13. If $n=2$, Theorem 4.11 can be strengthened. For $n=2$, the following three assertions are equivalent.

(a) Strategy $k$ is a strict Nash equilibrium of $\tilde{A}$.
(b) The point $\mathbf{e}_k$ is stochastically asymptotically stable.
(c) For some $\mathbf{x} \in \text{int}\,\Delta$, $P_\mathbf{x}\{\lim_{t\to\infty} X(t) = \mathbf{e}_k\} > 0$.

This is in contrast to the deterministic case where a pure strategy can be asymptotically stable without being a strict equilibrium.

For the proof, it remains to show that (c) implies (a). Suppose $P_\mathbf{x}\{X(t) \to \mathbf{e}_1\} > 0$. By Theorem 4.11(b), $\tilde{a}_{11} \geq \tilde{a}_{21}$, and if $\tilde{a}_{11} = \tilde{a}_{21}$, then $\tilde{a}_{22} < \tilde{a}_{12}$. But if $\tilde{a}_{11} = \tilde{a}_{21}$ and $\tilde{a}_{22} < \tilde{a}_{12}$, then $\{X(t)\}$ is null recurrent; see Example 4.3. This contradicts the assumption that $P_\mathbf{x}\{X(t) \to \mathbf{e}_1\} > 0$, and it follows that $\tilde{a}_{11} > \tilde{a}_{21}$.

COROLLARY 4.14. *If $n \in \{2,3\}$ and $A$ is conditionally positive definite, then $\{X(t)\}$ is transient.*

PROOF FOR $n=3$. If none of the columns of $\tilde{A}$ is proportional to $\mathbf{1}$, then transience follows from Corollary 4.8. Suppose now that, for example, $\tilde{a}_{11} = \tilde{a}_{21} = \tilde{a}_{31}$. Since $\tilde{A}$ is conditionally positive definite, $0 < \tilde{a}_{ii} - \tilde{a}_{ji} - \tilde{a}_{ij} + \tilde{a}_{jj}$ for $i \neq j$. Hence $\tilde{a}_{12} < \tilde{a}_{22}$ and $\tilde{a}_{13} < \tilde{a}_{33}$. Moreover, $0 < \tilde{a}_{22} - \tilde{a}_{32} - \tilde{a}_{23} + \tilde{a}_{33}$, so that $\tilde{a}_{22} > \tilde{a}_{32}$ or $\tilde{a}_{33} > \tilde{a}_{23}$. It follows that strategy 2 or strategy 3 is a strict Nash equilibrium of $\tilde{A}$. Thus, by Theorem 4.11, $\mathbf{e}_2$ or $\mathbf{e}_3$ is stochastically asymptotically stable and $\{X(t)\}$ is transient. $\square$

EXAMPLE 4.15. Consider the rock–scissors–paper game with payoff matrix

$$A = \begin{pmatrix} 0 & -a_1 & a_2 \\ a_2 & 0 & -a_1 \\ -a_1 & a_2 & 0 \end{pmatrix}, \qquad a_1, a_2 > 0.$$

Thus strategy 3 beats strategy 2, 2 beats 1 and 1 beats 3. If $a_1 > a_2$, then for all $\mathbf{y} \in \mathbb{R}^3 \setminus \{\mathbf{0}\}$ with $y_1 + y_2 + y_3 = 0$,

$$\mathbf{y}^T A \mathbf{y} = \frac{a_1 - a_2}{2} \{y_1^2 + y_2^2 + y_3^2\} > 0.$$

Hence, according to Corollary 4.14, $\{X(t)\}$ is transient. On the other hand, if $a_1 < a_2$, $\{X(t)\}$ is positive recurrent, provided $\sigma_1, \sigma_2, \sigma_3$ are sufficiently small; see Imhof (2005b). Moreover, if $a_1 < a_2$ and $\sigma_1 = \sigma_2 = \sigma_3 > 0$, Corollary 3.8 yields that the invariant distribution is the Dirichlet distribution with parameter $(a_2 - a_1)/(3\sigma_1^2)\mathbf{1}$. Finally, if $a_1 = a_2$ and $\sigma_1 = \sigma_2 = \sigma_3 > 0$, then Corollary 3.5 and 3.10 each imply that $\{X(t)\}$ is not positive recurrent. To see that this follows from Corollary 3.5, note that the stochastic replicator dynamics remains unchanged if a constant is added to every entry in a column of $A$. In view of Corollary 3.12 we conjecture that $\{X(t)\}$ is null recurrent in this case.

COROLLARY 4.16. *Suppose $A$ satisfies (3.3) with some $\gamma \in \mathbb{R}$. Then each of the following conditions is sufficient for transience of $\{X(t)\}$.*

(i) $\gamma < 0$ and $n \in \{2, 3\}$.
(ii) *There exists* $\mathbf{p} \in \Delta \setminus \{\mathbf{e}_1, \ldots, \mathbf{e}_n\}$ *with* $\{\tilde{A}\mathbf{p}\}_1 = \cdots = \{\tilde{A}\mathbf{p}\}_1$ *and* $\gamma < 0$.
(iii) *There does not exist* $\mathbf{p} \in \Delta$ *with* $\{\tilde{A}\mathbf{p}\}_1 = \cdots = \{\tilde{A}\mathbf{p}\}_1$.

PROOF. This follows from Remark 3.7, Corollary 4.14 and Theorems 4.1 and 4.5. □

**5. Related results.** The aim of this section is to put our results into perspective by discussing the connection with corresponding results for the deterministic replicator dynamics. We concentrate on the following fundamental properties of the deterministic model (1.2).

(A) Every Nash equilibrium of $A$ is a rest point.
(B) Every strict Nash equilibrium of $A$ is asymptotically stable.
(C) Every stable rest point is a Nash equilibrium of $A$.
(D) If an interior orbit converges, its limit is a Nash equilibrium of $A$.
(E) Every evolutionarily stable strategy (ESS) of $A$ is asymptotically stable.
(F) If (1.2) is permanent, that is, if there exists a compact set $K \subset \mathrm{int}\,\Delta$ such that for every initial state $\boldsymbol{\xi}(0) = \mathbf{x} \in \mathrm{int}\,\Delta$, $\boldsymbol{\xi}(t) \in K$ for $t$ large enough, then there exists a unique interior rest point and for any initial state in $\mathrm{int}\,\Delta$, the time averages $T^{-1} \int_0^T \boldsymbol{\xi}(t)\, dt$ converge to that rest point.

Parts (A)–(D) constitute the folk theorem of evolutionary game theory. For proofs of all of these statements, see Hofbauer and Sigmund (1998). The stochastic replicator dynamics (1.4) has no absorbing points except for the vertices $\mathbf{e}_1, \ldots, \mathbf{e}_n$, so that part (A) is in marked contrast to the stochastic situation. Note that the converse of (A) holds neither for the deterministic nor for the stochastic replicator dynamics.

Theorem 4.11(a) shows that part (B) essentially carries over to the stochastic case if the payoff matrix $A$ is replaced by $\tilde{A}$. In contrast to what might be expected from the deterministic case, the assumption that a given strategy is a strict Nash equilibrium does not imply almost sure convergence to the equilibrium for suitable initial values. Indeed, if $\tilde{A}$ describes a coordination game with two strategies, then for every initial state $\mathbf{x} \in \text{int}\, \Delta$, $P_{\mathbf{x}}\{X(t) \to \mathbf{e}_1\} < 1$. Theorem 4.11(b) can be regarded as a weak analog of part (C). A sufficient condition for $P_{\mathbf{x}}\{X(t) \to \mathbf{e}_k\} = 1$ for all $\mathbf{x} \in \text{int}\, \Delta$ is given in Imhof (2008).

For the process $\{X^{(s)}(t)\}$ determined by the Stratonovich stochastic differential equation (1.6), we obtain from Theorem 4.11(a) (applied to the corresponding Itô equation) that $\mathbf{e}_k$ is stochastically asymptotically stable if $a_{kk} > \max_{j \neq k} a_{jk}$. That the inequality

$$a_{kk} > \max_{j \neq k} a_{jk} + \frac{1}{8}\left(\sum_{j \neq k} \sigma_j^2 - \frac{(n-3)^2}{\sum_{j \neq k} \sigma_j^{-2}}\right)$$

is sufficient for stochastic asymptotic stability of $\mathbf{e}_k$ has been shown by Khasminskii and Potsepun (2006). In view of Example 4.3, the sufficient condition $a_{kk} > \max_{j \neq k} a_{jk}$ cannot be weakened to $a_{kk} \geq \max_{j \neq k} a_{jk}$.

Theorem 4.11(b) also yields the following analog of part (D). If for some $\mathbf{x} \in \text{int}\, \Delta$ and some $\mathbf{p} \in \Delta$, $P_{\mathbf{x}}\{X(t) \to \mathbf{p}\} > 0$, then, by Theorem 4.3 of Imhof (2005a), $\mathbf{p} \in \{\mathbf{e}_1, \ldots, \mathbf{e}_n\}$, so that, by Theorem 4.11(b), $\mathbf{p}$ is a Nash equilibrium of $\tilde{A}$. However, the natural stochastic analog of (D) is provided by Theorem 2.1(a).

Concerning part (E) note that in Example 4.3, strategy $\mathbf{e}_1$ is an ESS of $\tilde{A}$, which is not stochastically asymptotically stable. However, the existence of an interior ESS of $A$ implies positive recurrence and that the invariant distribution puts most mass near the ESS, provided that $\sigma_1, \ldots, \sigma_n$ are small enough; see Imhof (2005a). Since the trajectories do not converge to the ESS nor to any other point, it is natural to study convergence of their time averages. By Theorem 3.1(a), the averages converge a.s. to the unique interior Nash equilibrium of $\tilde{A}$.

Theorem 3.1(a) can also be seen as showing that an analog of the conclusion of part (F) holds when $\{X(t)\}$ is positive recurrent. The connection between permanence and recurrence is examined in Benaïm, Hofbauer and Sandholm (2008).

**Acknowledgment.** We are grateful to the referees for their constructive comments.

Department of Mathematics  
University of Vienna  
A-1090 Vienna  
Austria  

Department of Statistics  
and Hausdorff Center for Mathematics  
Bonn University  
53113 Bonn  
Germany  
E-mail: limhof@uni-bonn.de